\documentclass[a4paper,12pt]{article}
\usepackage{a4wide}
\usepackage{theorem}
\usepackage{amsmath}
\usepackage{array}
\usepackage{amssymb}
\usepackage{amsfonts}
\usepackage[english]{babel}
\usepackage{epsf}
\usepackage{epsfig}
\usepackage{graphicx}

\newtheorem{theo}{\indent Theorem \newline}[section]

\newtheorem{prop}[theo]{\indent Proposition\newline}

 \def\N{{\mathbb{N}}}
\def\Z{{\mathbb{Z}}}

\def\R{{\mathbb{R}}}
\def\C{{\mathbb{C}}}

\def\E{{\mathbb{E}}}

\usepackage{eufrak}

\newlength{\indentation}%
\makeatletter
\newcommand\@makefntextsans[1]{%
    \parindent 0em%
    \noindent%
    \hb@xt@0em{\hss}%
    #1}
\def\footnotetextsans{%
     \@ifnextchar [\@xfootnotenextsans%
       {\@footnotetextsans}}
\def\@xfootnotenextsans[#1]{%
  \begingroup%
     \csname c@\@mpfn\endcsname #1\relax%
  \endgroup%
  \@footnotetextsans}
\long\def\@footnotetextsans#1{\insert\footins{%
    \reset@font\footnotesize%
    \interlinepenalty\interfootnotelinepenalty%
    \splittopskip\footnotesep%
    \splitmaxdepth \dp\strutbox \floatingpenalty \@MM%
    \hsize\columnwidth \@parboxrestore%
    \color@begingroup%
      \@makefntextsans{%
        \rule\z@\footnotesep\ignorespaces#1\@finalstrut\strutbox}
    \color@endgroup}}
\makeatother

\begin{document}

\title{Topology of random real hypersurfaces}
\author{Jean-Yves Welschinger}
\maketitle

\makeatletter\renewcommand{\@makefnmark}{}\makeatother
\footnotetextsans{Keywords: random polynomials, real algebraic maifolds, random matrices.}
\footnotetextsans{AMS Classification : 14P25, 60D05.
}

{\bf Abstract:}

These are notes of the mini-course I gave during the CIMPA summer school at Villa de Leyva, Colombia, in July $2014$.
The subject was my joint work with Damien Gayet on the topology of random real hypersurfaces, restricting myself to
the case of projective spaces and focusing on our lower estimates. Namely, we estimate from (above and) below the mathematical
expectation of all Betti numbers of degree $d$ random real projective hypersurfaces. For any closed connected hypersurface $\Sigma$
of $\R^n$, we actually estimate from below the mathematical expectation of the number of connected components of these 
degree $d$ random real projective hypersurfaces which are diffeomorphic to $\Sigma$. 

\section{Random real polynomials}

\subsection{In one variable}

Let $P \in \R_d [X]$ be a polynomial in one variable, of degree $d$ and with real coefficients. Let $V_P = \{ x \in \R \, \vert \, P(x) = 0 \}$
be the  set of its real roots. Everybody knows the following
\begin{theo}
\label{theofond}
Let $P \in \R_d [X]$, then $0 \leq \# V_P \leq d$ and moreover $\# V_P \equiv d \mod (2)$ provided $P$ is generic enough. 
\end{theo}

This Theorem \ref{theofond} already raises a question which is going to be the main theme of this course. \\

{\bf Question $1$:}  What is the typical number of roots of $P$, choosing $P$ at random?\\

The mathematical expectation of this number of real roots reads as the average

$$\E (\# V_P) = \int_{\R_d [X]} (\# V_P) d\mu (P),$$
where $\mu$ denotes some probability measure on $\R_d [X]$.\\

{\bf First answer:}

A first answer to Question $1$ has been given by M. Kac in the $40's$

\begin{theo}[M. Kac, $1943$, \cite{Kac}]
\label{theoKac}
$$\E (\# V_P) \sim_{d \to + \infty} \frac{2}{\pi} \log(d).$$
\end{theo}
In order to provide this answer, Kac did consider that the space $\R_d [X]$ of polynomials is Euclidean, a canonical orthonormal basis
being given by the monomials $1, X , X^2, \dots, X^d$. Now, since this space is Euclidean, it carries a canonical 
probability measure, the Gaussian measure associated to its scalar product. The latter reads
$$d\mu (P) = \frac{1}{\sqrt{\pi}^{d+1}} \exp (-\| P \|^2) dP,$$
where $d+1$ corresponds to the dimension of $\R_d [X]$ and $dP$ to its Lebesgue measure which is associated
to the scalar product but has infinite volume. This Gaussian measure is thus the Lebesgue measure weighted with some exponential
which reduces its total volume to one. It has the great properties to be a product measure which is invariant under
the orthogonal group. \\

{\bf Second answer:}

A second answer to Question $1$ has been given by E. Kostlan. 

\begin{theo}[Kostlan, Shub-Smale $1993$, \cite{Ko}, \cite{SS}]
\label{theoKo}
For every $d>0$, 
$$ \E (\# V_P) = \sqrt{d}.$$
\end{theo}

In order to provide this answer, E. Kostlan also did equip the space $\R_d [X]$ of polynomials with some Gaussian measure, 
but associated to a different scalar product. For this new scalar product, an orthonormal basis is given by the monomials
$\sqrt{d \choose k} X^k$, $0 \leq k \leq d$. This scalar product turns out to be more natural geometrically. I will give a geometric
definition in \S \ref{subsecprob}, but let me already point out a nice property. 

The space $\R_d [X]$ is well known to be isomorphic to the space $\R_d^\text{hom} [X, Y]$ of homogeneous polynomials of degree $d$ in two variables and real coefficients. This isomorphism reads $X^k \in \R_d [X] \mapsto X^k Y^{d-k} \in \R_d^\text{hom} [X, Y]$. If we push forward
the new scalar product under this isomorphism, then we get one which is invariant under the action of the orthogonal group of the plane, 
by composition on the right. That is, for every $Q \in \R_d^\text{hom} [X, Y]$ and every $h \in O_2 (\R)$, $\| Q \| = \| Q \circ h^{-1} \|$. 

Let me explain the proof of Kostlan, which also recovers the result of Kac. 

{\bf Proof: (see \cite{EdKo})}

Let us fix the isomorphism $(a_0 , \dots , a_d) \in \R^{d+1} \mapsto \sum_{i=0}^d a_i X^i \in \R_d [X]$
and focus on two objects of 
$\R^{d+1} $. First, the unit sphere $S^{d} $ and for every $\underline{a} = (a_0 , \dots , a_d) \in  S^{d} $, let us denote by
$\lambda_{\underline{a} }$ the linear form $(y_0 , \dots , y_d) \in \R^{d+1} \mapsto \sum_{i=0}^d a_i y_i \in \R$. 

Secondly, let us consider the curve $\tilde{\gamma} : t \in \R \mapsto (1 , t , \dots , t^d) \in \R^{d+1}$ in the case of Kac, or
$\tilde{\gamma} : t \in \R \mapsto (1 , \sqrt{d \choose 1} t , \dots ,\sqrt{d \choose d}  t^d) \in \R^{d+1}$ in the case of Kostlan.

If $P = \sum_{i=0}^d a_i X^i \in \R_d [X]$, then, as a function on the real line, $P = \lambda_{\underline{a} } \circ \tilde{\gamma} $ (in the case
of Kac), so that $V_P \cong \ker \lambda_{\underline{a} } \cap \text{Im} (\tilde{\gamma} )$.


The observation of Kostlan is then the following.
\begin{theo}[\cite{EdKo}]
\label{theoEdKo}
$$\E (\# V_P) = \frac{1}{\pi} \text{length} (\gamma), \text{ where } \gamma : t \in \R \mapsto  \frac{\tilde{\gamma} (t)}{\| \tilde{\gamma} (t) \|} \in S^{d}.$$
\end{theo}

{\bf Proof : (see \cite{EdKo})}

This is Crofton's formula, the length of a curve is the average of the number of intersection points with the hyperplanes. Note that the formula
is obvious when the curve is a closed geodesic on the sphere. This formula follows from the fact that it also holds true for a piece of such a geodesic,
since every smooth curve can then be approximated by some piecewise geodesic curve. $\square$

The end of the proof of Theorem \ref{theoEdKo} is just a computation of the length of the curve $\gamma$, which gives
$\text{length} (\gamma) \sim_{d \to + \infty} 2 \log(d)$ in the case of Kac (a bit tough) and
$\text{length} (\gamma) = \pi \sqrt{d}$ in the case of Kostlan (easy). $\square$

\subsection{In several variables}

What about polynomials in several variables?

If $P \in \R_d [X_1 , \dots , X_n]$ is a polynomial in $n$ variables, degree $d$ and real coeffcients, then
 $V_P = \{ x \in \R^n \, \vert \, P(x) = 0 \}$ is no more a finite set, but rather an affine real algebraic hypersurface. 
It is not compact in general, but has a standard compactification. Namely, this space of polynomial is again canonically isomorphic
to the space $\R_d^{\text{hom}} [X_0 , \dots , X_n]$ of homogeneous polynomials of degree $d$, $n+1$ variables and real coefficients. 
This isomorphism reads $X_1^{\alpha_1}\dots X_n^{\alpha_n} \in \R_d [X_1 , \dots , X_n] \mapsto X_0^{d- \alpha_1 - \dots - \alpha_n}
X_1^{\alpha_1}  \dots X_n^{\alpha_n} \in \R_d^{\text{hom}} [X_0 , \dots , X_n]$ and if $Q \in \R_d^{\text{hom}} [X_0 , \dots , X_n]
\setminus \{ 0 \}$, then 
$V_Q = \{ x \in \R P^n \, \vert \, Q(x) = 0 \}$ is a compact hypersurface, smooth for generic polynomials and which then contains $V_P$ as
a dense subset. I will come back to projective spaces in \S \ref{subsecproj}. 

Again, the topology of $V_Q$ depends on the choice of $Q$, as in one variable. 
For example, in degree $d=2$ and $n=3$ variables, $V_Q$ is a quadric surface which may be empty, homeomorphic to a sphere in the case
of the ellipsoid or to a torus
in the case of the hyperboloid. If we denote, for every $i \in \{0 , \dots , n-1 \}$,
by $b_i (V_Q ; \Z/2\Z)= \dim H_i  (V_Q ; \Z/2\Z)$ the $i$-th Betti number with $\Z/2\Z$ coefficients of $V_Q$, then
\begin{theo}[Smith-Thom's inequality, $1965$, \cite{Thom}]
\label{theoSmTh}
$$0 \leq \sum_{i=0}^{n-1} b_i (V_Q ; \Z/2\Z) \leq \sum_{i=0}^{2n-2} b_i (V_Q^{\C} ; \Z/2\Z) = d^n + o(d^n).$$
Moreover, $\sum_{i=0}^{n-1} b_i (V_Q ; \Z/2\Z) \equiv \sum_{i=0}^{2n-2} b_i (V_Q^{\C} ; \Z/2\Z) \mod(2)$, provided 
$Q \in \R_d^{\text{hom}} [X_0 , \dots , X_n]$ is generic enough for $V_Q$ to be smooth. 
\end{theo}
In Theorem \ref{theoSmTh}, $V_Q^{\C} = \{ x \in \C P^n \, \vert \, Q(x) = 0 \}$ denotes the set of complex roots of $Q$ in the complex projective space.
This Theorem \ref{theoSmTh} extends Theorem \ref{theofond}, which corresponds to the case $n=1$. The case $n=2$ was also previously
known as the (famous in real algebraic geometry) Harnack-Klein's inequality, see \cite{Har}, \cite{Klein}.

Again, this raises the question\\

{\bf Question $2$:} What is the typical topology of $V_Q$, choosing $Q$ at random in $\R_d^{\text{hom}} [X_0 , \dots , X_n]$?

e.g. which Betti numbers to expect?\\

Let me give a formulation of our joint results with Damien Gayet. 

For every $i \in \{0 , \dots , n-1 \}$, let me set
\begin{eqnarray*}
b_i (V_Q ; \R) & = & \dim H_i  (V_Q ; \R),\\
m_i (V_Q) & = & \inf_{f \text{ Morse on } V_Q} \# \text{Crit}_i (f).
\end{eqnarray*}
Here, $\text{Crit}_i (f)$ denotes the number of critical points of index $i$ of the Morse function $f$. Recall that a real function
of class $C^2$ is said to be Morse if and only if all of its critical points are non-degenerated. This means that the Hessian of this function
at all of its critical points is a non-degenerated quadratic form. The index of such a quadratic form is then the maximal dimension of a linear
subspace of the tangent space at the critical point where it restricts to a negative definite one, see \cite{Milnor}. The latter is called the index
 of the critical point. 

It follows from Morse theory that $ b_i (V_Q ; \R) \leq b_i (V_Q ; \Z/2\Z) \leq m_i (V_Q)$, see \cite{Milnor}. 

The mathematical expectations for these Betti or Morse numbers read as the averages
\begin{eqnarray*}
\E (b_i) & = & \int_{\R_d^{\text{hom}} [X_0 , \dots , X_n]} b_i (V_Q ; \R)  d\mu (Q),\\
\E (m_i) & = & \int_{\R_d^{\text{hom}} [X_0 , \dots , X_n]} m_i (V_Q)  d\mu (Q).
\end{eqnarray*}

The probability measure $\mu$ we consider extends the one considered in Theorem \ref{theoKo}. It is the Gaussian measure associated to 
the scalar product for which the monomials $\sqrt{\frac{(d+n) !}{n! \alpha_0 ! \dots \alpha_n !}} X_0^{\alpha_0} \dots X_n^{\alpha_n} $,
$\alpha_0 + \dots + \alpha_n = d$, define an orthonormal basis. Again, the action of the orthogonal group of the $(n+1)$-dimensional Euclidean space by composition on the right preserves this scalar product. That is, for every $Q \in \R_d^\text{hom} [X_0 , \dots , X_n]$ and every $h \in O_{n+1} (\R)$, $\| Q \| = \| Q \circ h^{-1} \|$. Let me finally observe that the coefficient $(d+n) !$ instead of $d!$ in the numerator of the mononials has only the effect to rescale the scalar product and does not affect the results. We will see in 
\S \ref{subsecprob} how this scalar product shows up.

\begin{theo}[joint with Damien Gayet, \cite{GaWe3}, \cite{GaWe4}]
\label{theoGaWe}
There exist (universal) constants $c_i^+ , c_i^-$ such that
$$c_i^- \leq \liminf_{d \to + \infty} \frac{\E (b_i)}{\sqrt{d}^n \text{Vol}_{\text{FS}} \R P^n} \leq \limsup_{d \to + \infty} \frac{\E (m_i)}{\sqrt{d}^n \text{Vol}_{\text{FS}} \R P^n} \leq c_i^+.$$
\end{theo}

Unfortunately, by lack of time, I will only explain the proof of the lower estimates given by Theorem \ref{theoGaWe} in this course. The term
$\text{Vol}_{\text{FS}} \R P^n$ denotes the total volume of the real projective space for the Fubini-Study metric, see \S \ref{subsecprob}.
Though it is some constant, I distinguish it from $c_i^\pm$. In fact, in \cite{GaWe3}, \cite{GaWe4}, we not only prove Theorem \ref{theoGaWe}
for projective spaces, but for any smooth real projective manifold. The term $\text{Vol}_{\text{FS}} \R P^n$ has then to be replaced by the total
K\"ahlerian volume of the real locus of the manifold, for the K\"ahler metric induced by the curvature form of a metric with positive curvature
chosen on some ample real line bundle, the tensor powers of which we consider random sections. The constants $c_i^+ , c_i^-$ are, they, unchanged and only depend on $i$ and $n$, see \S \ref{subseccsts}.

\subsection{The universal constants $c_i^+ , c_i^-$ }
\label{subseccsts}

Let me tell you more about these universal constants $c_i^+ , c_i^-$. 

The constant $c_i^+ $ is related to random symmetric matrices. Namely $c_i^+ = \frac{1}{\sqrt{\pi}} e_\R (i , n-1-i)$, where
\begin{eqnarray}
\label{intsym}
e_\R (i , n-1-i) & = &\int_{\text{Sym}(i , n-1-i ; \R)} \vert \det (A) \vert d\mu (A).
\end{eqnarray}
Here, $\text{Sym}(i , n-1-i ; \R)$ denotes the open cone of non-degenerated symmetric matrices of size $(n-1) \times (n-1)$, signature
$(i , n-1-i)$ and real coefficients. It is included in the vector space $\text{Sym}(n-1 ; \R)$ of real symmetric matrices of size $(n-1) \times (n-1)$.
The latter is Euclidean, equipped with the scalar product $(A,B) \in \text{Sym}(n-1 ; \R)^2 \mapsto \frac{1}{2} tr(AB) \in \R$, see \cite{Mehta}. 
So again this space inherits some Gaussian measure $\mu$, which is the one we consider in the integral (\ref{intsym}).

In particular, 
$$ \sum_{i=0}^{n-1} c_i^+ = \frac{1}{\sqrt{\pi}} \E (\vert \det (A) \vert ) = \frac{1}{\sqrt{\pi}} \int_{\text{Sym}(n-1; \R)} \vert \det (A) \vert d\mu (A).$$

\begin{theo}[joint with Damien Gayet, \cite{GaWe3}]
\label{theoGaWesym}
\begin{enumerate}
\item $ \sum_{i=0}^{n-1} c_i^+\sim_{n \to + \infty} \frac{2\sqrt{2}}{\pi} \Gamma(\frac{n+1}{2})$, where $ \Gamma$ denotes Euler's function. 
\item For every $\alpha \in [0 , \frac{1}{2} [$, there exists $c_\alpha >0$ such that for $n$ large enough 
$\sum_{i=0}^{\lfloor \alpha n \rfloor} c_i^+ \leq \exp (-c_\alpha n^2)$. $\square$
\end{enumerate}
\end{theo}
Note that the first part of Theorem \ref{theoGaWesym} was known for $n$ even, see \cite{Mehta} and \cite{GaWe3}, while the second part
quickly follows from some large deviation estimates established in \cite{BAG}. 

For instance, for $i=0$, $b_i (V_Q ; \R)  = b_0 (V_Q ; \R) $ denotes the number of connected components of $V_Q$, Theorem \ref{theoGaWe}
estimates the expected number of connected components of $V_Q$ and $c_i^+ = c_0^+ $ provides the upper estimate. This constant
more than exponentially decreases as the dimension $n$ grows to $+ \infty$. 

As for the constant $c_i^-$, we set

$${\cal H}_n = \{ \text{closed connected hypersurfaces of } \R^n \} /\text{diffeomorphisms}.$$
For every $[\Sigma] \in {\cal H}_n $, we set $b_i (\Sigma ; \R) = \dim H_i (\Sigma ; \R) $ and associate some positive constant $c_{[\Sigma] }$,
see \S \ref{subsecproof}. This constant $c_{[\Sigma] }$ is defined via some quantitative transversality, but turns out at the end to bound from below the expected number of connected components of $V_Q$ that are diffeomorphic to $\Sigma$, which is what we actually estimate from below, see \S \ref{subsecproof}. Then, $c_i^-$ is the infinite serie
$$c_i^- = \sum_{[\Sigma] \in {\cal H}_n} c_{[\Sigma] } b_i (\Sigma ; \R).$$
This serie converges since it is bounded from above by $c_i^+$.

\begin{prop}[joint with Damien Gayet, \cite{GaWe4}]
\label{propGaWe}
For every $i \in \{ 0 , \dots , n-1 \}$, $c_i^- \geq c_{[S^i \times S^{n-1-i}]} \geq \exp (-\exp{(5n + 69)})$. $\square$
\end{prop}
 Indeed, the product of the $i$-dimensional unit sphere with the $(n-1-i)$-dimensional unit sphere turns out to embed as a closed connected
hypersurface of $\R^n$. The $i$-th Betti number of this hypersurface is one and we will see in \S \ref{subsecproof} that the constant 
$c_{[\Sigma] }$ is actually explicit, so that it can be estimated for this product of spheres, see \cite{GaWe4}.

\section{The lower estimates}

\subsection{Projective spaces}
\label{subsecproj}

Recall that the $n$-dimensional projective space is by definition the space of one-dimensional linear subspaces of the affine
$(n+1)$-dimensional space. That is, 
\begin{eqnarray*}
\C P^n &=& \{ \text{space of lines in } \C^{n+1} \} \\
 &=& \C^{n+1} \setminus \{ 0 \} /x \sim \lambda x, \forall \lambda \in \C^*,
\end{eqnarray*}
and likewise,
\begin{eqnarray*}
\R P^n &=& \{ \text{space of lines in } \R^{n+1} \} \\
 &=& \R^{n+1} \setminus \{ 0 \} /x \sim \lambda x, \forall \lambda \in \R^*.
\end{eqnarray*}
The points in $\C P^n$ are represented by their homogeneous coordinates $[x_0 : \dots : x_n]$, where $x_0 , \dots , x_n \in \C$ do not all
vanish, being understood that for every $\lambda \in \C^*$, $[x_0 : \dots : x_n] = [\lambda x_0 : \dots : \lambda x_n]$.

These complex projective spaces are smooth compact complex manifolds without boundary. They are covered by $n+1$ standard affine charts. 
Namely, for every $i \in \{ 0 , \dots , n \}$, set $U_i = \{ [x_0 : \dots : x_n] \in \C P^n \, \vert \, x_i \neq 0 \}$. This dense open subset
$U_i$ corresponds to the lines of $\C^{n+1}$ that are not contained in the hyperplane $\{ (x_0 , \dots , x_n) \in \C^{n+1} \, \vert \, x_i = 0 \}$.
Every such line intersects the affine hyperplane $\{ (x_0 , \dots , x_n) \in \C^{n+1} \, \vert \, x_i = 1 \}$ at exactly one point, defining the chart
$$\phi_i :  [x_0 : \dots : x_n] \in U_i \mapsto (\frac{x_0}{x_i}, \dots, \frac{x_{i-1}}{x_i}, \frac{x_{i+1}}{x_i}, \dots,\frac{x_n}{x_i}) \in \C^n.$$

\subsection{Line bundles}
\label{subsecline}

The projective space $\C P^n$ is the space of lines of $\C^{n+1}$, so that every point $x= [x_0 : \dots : x_n]$ of $\C P^n$ represents
a complex line $\gamma_x \subset \C^{n+1}$, the line generated by $(x_0 , \dots , x_n)$ in $\C^{n+1}$. The collection of all lines
$\gamma_x$, $x \in \C P^n$, defines what is called a holomorphic line bundle $\gamma$ over $\C P^n$. 
It is in particular a complex manifold equipped with a holomorphic submersion onto the base $\C P^n$, see \cite{Dem}, \cite{GH}. 
Since all these lines are included
in $\C^{n+1}$, the tautological line bundle $\gamma$ is a subline bundle of the trivial vector bundle $\C P^n \times \C^{n+1} \to \C P^n$
of rank $n+1$.

Now, every vector space comes with its dual space, the space of linear forms over it. This defines the dual bundle
$\gamma^* = \{ \text{linear forms on } \gamma \} \to \C P^n$. Likewise, for every $d >0$, I denote by $\gamma_d^*$ the space
of homogeneous forms of degree $d$ on $\gamma$, so that $\gamma_1^* = \gamma^*$. Again, all these define holomophic line bundles over $\C P^n$.
Note that another standard notation for these bundle is $\gamma = {\cal O}_{\C P^n} (-1)$, $\gamma_d^* = {\cal O}_{\C P^n} (d)$.

We denote by $H^0 (\C P^n ; \gamma_d^*)$ the space of global holomorphic sections of the bundle $\gamma_d^*$, that is the space
of holomorphic maps $s : \C P^n \to  \gamma_d^*$ such that $\pi \circ s = \text{id}_{\C P^n}$, where $\pi : \gamma_d^* \to \C P^n$
denotes the tautological projection. Hence, for every point $x \in \C P^n$, $s(x)$ denotes a homogeneous form of degree $d$ on 
the complex line $\gamma_x$.

Now, complex homogeneous polynomials of degree $d$ in $n+1$ variables define homogeneous functions of degree $d$ on $\C^{n+1}$.
These thus restrict to homogeneous functions of degree $d$ on every line  $\gamma_x$, whatever $x \in \C P^n$ is. As a consequence, 
these  complex homogeneous polynomials of degree $d$ define global holomorphic sections of the bundle $\gamma_d^*$ so that we
get an injective morphism $\C_d^{\text{hom}} [X_0 , \dots , X_n] \hookrightarrow H^0 (\C P^n ; \gamma_d^*)$. It is not that hard to prove that
this injective morphism is also surjective, but requires though two theorems in complex analysis, namely Hartog's theorem and the decomposition of entire functions into power series, see \cite{Dem}.

{ \bf Upshot:} It is important here to understand that a homogeneous polynomial $Q \in \C_d^{\text{hom}} [X_0 , \dots , X_n]$ does not
define a holomorphic function $\C P^n \to \C$ (any such function would be constant due to maximum's principle). Its vanishing subset
in $\C^{n+1} \setminus \{ 0 \}$ is a cone, and thus defines on the quotient $\C P^n$ the hypersurface 
$V_Q^{\C} = \{ x \in \C P^n \, \vert \, Q(x) = 0 \}$, provided $Q$ is of positive degree. But the other level sets of $Q$ in 
$\C^{n+1} \setminus \{ 0 \}$ are not left invariants under homotheties and thus do not pass to the quotient $\C P^n$.

What is true is that these polynomials $Q \in \C_d^{\text{hom}} [X_0 , \dots , X_n]$ define global sections of $\gamma_d^*$, and $V_Q^{\C}$ coincides with the vanishing locus of these as sections of $\gamma_d^*$.

\subsection{Fubini-Study metric}
\label{subsecFS}

Let me now equip $\C^{n+1}$ with its standard Hermitean product, defined for every $v=(v_0 , \dots , v_n)$ and $w=(w_0 , \dots , w_n)$
in $\C^{n+1}$ by $h(v,w) = \sum_{i=0}^n v_i \overline{w}_i \in \C$. 

It restricts on every line $\gamma_x$ of $\C^{n+1}$ to a Hermitean product $h$. This is called a Hermitean metric on the line bundle 
 $\gamma$. It also induces then a Hermitean metric $h_d$ on all the line bundles $\gamma_d^*$, $d>0$. Indeed, if $x \in \C P^n$ and
$s(x) \in \gamma_d^*\vert_x$, then $s(x) : \gamma_x \to \C$ is a homogeneous form of degree $d$ and
we set 
$$\| s(x) \| = \frac{\vert s(x)(v) \vert}{\| v \|^d},$$
where this definition does not depend on the choice of $v \in \gamma_x \setminus \{ 0 \}$.\\

{\bf Fundamental example:} Let us compute the pointwise Fubini-study norm of $Q = X_0^d \in \R_d^{\text{hom}} [X_0 , \dots , X_n]$,
viewed as a section of $\gamma_d^*$.

I restrict myself to $U_0 \cong \C^n$, since it vanishes outside of $U_0$. Let $x = [1 : x_1 : \dots : x_n] \in U_0$. Then, 
$v = (1 , x_1 , \dots , x_n)$ generates
$\gamma_x$ and $\| v \|^2 = h(v,v) = 1 + \sum_{i=1}^n \vert x_i \vert^2$. Since $Q(v) = \vert Q(v) \vert = 1$, 
we get
$$h_d (Q,Q)\vert_x = \frac{1}{(1 + \| x \|^2 )^d} = \exp (-d\log(1 + \| x \|^2)).$$

This means that the norm of $Q$ at the origin $[1 : 0 : \dots : 0] \in U_0$ equals one, but at every other point it decays exponentially fastly
to zero as the degree grows to $+ \infty$. Since $\log h_d (Q,Q)\vert_x = -d \| x \|^2 + o(\| x \|^2)$ near
$x = 0$, we deduce that the Fubini-Study norm of $Q$ gets concentrated in a ball of radius $\frac{1}{\sqrt{d}}$ centered at the origin. 
Such a section defined by $Q$ is called a peak section, and the scale $\frac{1}{\sqrt{d}}$ is a fundamental scale in K\"ahler geometry. 
Peak sections exist over any projective or Stein manifolds, following the theory of L. H\"ormander, see \cite{Horm}, \cite{MaMa}.

Finally, it is possible to define sections of $\gamma_d^*$ which peak near any point $x \in \C P^n$. Indeed, the group $GL_{n+1} (\C)$
acts by linear automorphisms of $\C^{n+1}$ and the unitary group $U_{n+1} (\C)$ even by isometries. These actions are transitive on lines
of $\C^{n+1}$ and thus they induce actions on $ \C P^n$ which  are transitive on points. Moreover, these actions lift to actions on 
$\gamma$ and thus on any line bundle $\gamma_d^*$, $d>0$.

For every $x \in \C P^n$, there exists $r \in U_{n+1} (\C)$ such that $x = r ([1 : 0 : \dots : 0] )$. Then, 
$Q \circ r^{-1} \in \C_d^{\text{hom}} [X_0 , \dots , X_n]$ defines a section of $\gamma_d^*$ which peaks near $ x $. 

\subsection{Implementation of affine hypersurfaces}
\label{subsechyp}

Let $\Sigma \subset \R^n$ be a closed hypersurface, not necessarily connected.
It is a theorem of H. Seifert, see \cite{Seif} or also \cite{Nash}, that there exists a polynomial $P \in \R [X_1 , \dots , X_n]$ of some degree $k$
such that $V_P = P^{-1} (0)$ contains a union of connected components $\widetilde{\Sigma}$ which is isotopic to $\Sigma$. This means that
there exists a path $(\phi_t)_{t \in [0,1]}$ of diffeomorphisms of $\R^n$ such that $\phi_0$ is the identity and 
$\phi_1 (\widetilde{\Sigma}) = \Sigma$. Note that this theorem of Seifert is similar to Stone-Weierstrass theorem, except that one needs
some approximation in $C^1$-norm. Note also that I could have immediately taken any polynomial $P$ for which zero is a regular value
and then defined $\Sigma$ to be any union of closed connected components of $V_P$. From now on, let me fix $P$ and denote by $\Sigma$ 
such a union of closed connected components of $V_P$.

There exists $R>0$ such that $\Sigma$ is included in the ball $B(0,R) \subset \R^n$ of radius $R$. Let me replace, for every $d>0$, $P$ by 
the polynomial $P_d = P(\sqrt{d} .)$. It is still a polynomial of degree $k$, whose coefficients are $O(\sqrt{d}^k)$. Indeed, if
$P = \sum_{(\alpha_1 , \dots , \alpha_n) \in \N^n} a_{\alpha_1 , \dots , \alpha_n} X_1^{\alpha_1} \dots X_n^{\alpha_n}$, then
$P_d = \sum_{(\alpha_1 , \dots , \alpha_n) \in \N^n} a_{\alpha_1 , \dots , \alpha_n}\sqrt{d}^{\alpha_1 + \dots + \alpha_n} X_1^{\alpha_1} \dots X_n^{\alpha_n}$.

Under the isomorphism $ \R_k [X_1 , \dots , X_n] \mapsto  \R_k^{\text{hom}} [X_0 , \dots , X_n]$, $P_d$ is mapped to the polynomial
$Q_d = \sum_{(\alpha_1 , \dots , \alpha_n) \in \N^n} a_{\alpha_1 , \dots , \alpha_n}\sqrt{d}^{\alpha_1 + \dots + \alpha_n} X_0^{k-\alpha_1 - \dots - \alpha_n} X_1^{\alpha_1} \dots X_n^{\alpha_n}$. After multiplication by $X_0^{d-k}$, it provides a section $\sigma_P = Q_d X_0^{d-k} \in
\R_d^{\text{hom}} [X_0 , \dots , X_n] = \R H^0 (\C P^n ; \gamma_d^*)$ which vanishes in the ball $B([1 : 0 : \dots : 0] , R/\sqrt{d}) \subset \R U_0
\subset \R P^n$ centered at the origin $[1 : 0 : \dots : 0]$ and of radius $R/\sqrt{d}$. Moreover, $\sigma_P^{-1} (0) \cap B([1 : 0 : \dots : 0] , R/\sqrt{d}) $ contains a union of components $\widetilde{\Sigma}$ such that the pair $(B([1 : 0 : \dots : 0] , R/\sqrt{d}) , \widetilde{\Sigma})$
gets diffeomorphic to $(\R^n , \Sigma)$. In addition, the pointwise Fubini-Study norm $h_d (\sigma_P , \sigma_P) = \| Q_d \|^2
\| X_0^{d-k} \|^2$ decays exponentially outside the origin as $d$ grows to $+ \infty$. This is indeed the case for
$ X_0^{d-k}$ as we saw in the previous paragraph, while $Q_d $ has fixed degree and coefficients $O(\sqrt{d}^k)$. We deduce that the
Fubini-Study norm of $\sigma_P$ is concentrated in the ball $(B([1 : 0 : \dots : 0] , R/\sqrt{d}) $.

Finally, after composition on the right by some suitable $r \in O_{n+1} (\R) \subset U_{n+1} (\C)$, we get for every $x \in \R P^n$
a section $\sigma_P \in \R_d^{\text{hom}} [X_0 , \dots , X_n]$ such that $\sigma_P^{-1} (0) \cap B(x , R/\sqrt{d}) $
contains a union of components $\widetilde{\Sigma}$ for which the pair $(B(x , R/\sqrt{d}) , \widetilde{\Sigma})$
gets diffeomorphic to $(\R^n , \Sigma)$ and such that the Fubini-Study norm of $\sigma_P$ exponentially decreases outside of 
this ball $B(x , R/\sqrt{d}) $. Note that since the radius of this ball converges to zero, the Riemannian metric of $\R P^n$ for which we
take the ball does not matter. We will however introduce the Fubini-Study metric of $\C P^n$ in the next paragraph. 

\subsection{The probability measure $\mu$ revisited}
\label{subsecprob}

Recall that I did introduce the projective spaces in \S \ref{subsecproj} and their tautological line bundles in \S \ref{subsecline}. These are
line subbundles of some trivial vector bundle. Let $x \in \C P^n$ and $\gamma_x \subset \C^{n+1}$ be the line it represents. 
Let $y \in \gamma_x \setminus \{ 0 \}$ and $p : \C^{n+1} \setminus \{ 0 \} \to \C P^n$ be the canonical projection. 
Then, the differential map $d_y p : T_y (\C^{n+1} \setminus \{ 0 \} ) = \C^{n+1}  \to T_x \C P^n$ contains $\gamma_x$ in its kernel
and restricts to an isomorphism $\gamma_x^\perp \to T_x \C P^n$, where $\gamma_x^\perp$ stands for the orthogonal
of $\gamma_x$ with respect to the standard Hermitean product of $\C^{n+1}$, see \S \ref{subsecFS}.

This hyperplane $\gamma_x^\perp$ does not depend on the choice of $y \in \gamma_x \setminus \{ 0 \}$, but the isomorphism
$d_y p \vert_{\gamma_x^\perp}$ does. By the way, the quotient of the trivial vector bundle $\C P^n \times \C^{n+1}$ by the tautological 
bundle $\gamma$ is not isomorphic to the tangent bundle $T \C P^n$. 

{\bf Exercise:} Prove that the latter tangent bundle $T \C P^n$ is rather
isomorphic to the bundle of morphisms from $\gamma$ to the former quotient bundle (while the quotient bundle is isomorphic to the
space of morphisms from the trivial line bundle to itself).

Let us now choose $y$ of norm one, so that it lies in the intersection of the unit sphere with $\gamma_x$. This intersection is a circle, 
the orbit of the action of the unitary group $U_1 (\C)$ by homothety. The circle fibration $S^{2n+1} \to \C P^n$ this action produces
is called the Hopf fibration. Still, the isomorphism $d_y p \vert_{\gamma_x^\perp}$ depends on the choice of $y \in 
S^{2n+1} \cap \gamma_x$, but up to an isometry, so that if we push forward under $d_y p$ the Hermitean product of $\gamma_x^\perp$, 
induced by restriction of the ambient one of $\C^{n+1}$, we get a well defined Hermitean product on $T_x \C P^n$, which does not
depend on the choice of $y \in  S^{2n+1} \cap \gamma_x$. The collection of all these Hermitean products on all tangent spaces of $\C P^n$
defines a Hermitean metric on $\C P^n$ called the Fubini-Study metric. The action of $U_{n+1} (\C)$ on $\C P^n$ we already discussed provides
isometries for this metric. 

{\bf Remark:} A projective line for this Fubini-Study metric has total area $\pi$ (exercise). We actually rescale in \cite{GaWe3}, \cite{GaWe4} this metric by a factor $\frac{1}{\sqrt{\pi}}$ to normalize this area to one. This is quite natural from another point of view, since this Hermitean Fubini-Study metric, which is actually a K\"ahler metric, also originate from the curvature form of the canonical connection associated to the Fubini-Study metric of
$\gamma$ introduced in \S \ref{subsecFS}. Since the cohomology class of this form is the first Chern class of the line bundle $\gamma$,
the volume of a projective line gets one for this metric. This Fubini-Study metric restricts to a Riemannian metric on $\R P^n$ and the quantity
$\text{Vol}_{\text{FS}} \R P^n$ in Theorem \ref{theoGaWe} is the total volume of $\R P^n$ for this Riemannian Fubini-Study metric. 

Now, since the line bundles $\gamma_d^*$ are equipped with Hermitean metrics and their base $\C P^n$ with some volume form $dx$,
induced by the Fubini-Study metric, the spaces $H^0 (\C P^n ; \gamma_d^*) = \C_d^{\text{hom}} [X_0 , \dots , X_n]$ of global holomorphic sections of these bundles inherit some $L^2$-Hermitean products, namely
$$(Q_1, Q_2) \in H^0 (\C P^n ; \gamma_d^*)^2 \mapsto \int_{\C P^n} h_d (Q_1, Q_2) dx \in \C.$$
These $L^2$-Hermitean products restrict on the spaces $\R H^0 (\C P^n ; \gamma_d^*) = \R_d^{\text{hom}} [X_0 , \dots , X_n]$ of real holomorphic sections to the  $L^2$-scalar products
$$(Q_1, Q_2) \in \R H^0 (\C P^n ; \gamma_d^*)^2 \mapsto \int_{\C P^n} h_d (Q_1, Q_2) dx \in \R.$$
Finally, now that the space of real homogeneous polynomials $\R_d^{\text{hom}} [X_0 , \dots , X_n]$ is again Euclidean, it inherits some
Gaussian measure $d \mu (P) = \frac{1}{\sqrt{\pi}^{N_d}} \exp (- \| P \|^2) dP$, where $N_d$ denotes the dimension of 
$\R_d^{\text{hom}} [X_0 , \dots , X_n]$ and $dP$ the Lebesgue measure associated to this $L^2$-scalar product.

{\bf Exercise:} The monomials $X_0^{\alpha_0} \dots X_n^{\alpha_n} $ are orthogonal to each other and in fact the probability measure $\mu$
is the one considered in Theorem \ref{theoGaWe}, so that
$\sqrt{\frac{(d+n) !}{n! \alpha_0 ! \dots \alpha_n !}} X_0^{\alpha_0} \dots X_n^{\alpha_n} $ is an orthonormal basis (provided the Fubini-Study metric on $\C P^n$ is normalized so that its total volume is one ; it is $\pi^n / n!$ for the metric just defined).\\

{\bf Example:} What is the $L^2$-norm of the section $\sigma_P$ we did construct in the previous \S \ref{subsechyp}?

Recall that $\sigma_P = Q_d X_0^{d-k} $, so that
\begin{eqnarray*}
\| \sigma_P \|^2 &=& \int_{\C P^n} h_d (\sigma_P , \sigma_P) dx \\
&=& \int_{U_0} h_d (\sigma_P , \sigma_P) dx \text{ since } dx (\C P^n \setminus U_0) = 0\\
&=& \int_{\C^n} \frac{\vert P(\sqrt{d} x) \vert^2}{(1 + \| x \|^2)^d} dx \\
&\sim_{d \to + \infty}& \frac{1}{d^n} \int_{\C^n} \vert P(y) \vert^2 \exp (-\| y \|^2) dy 
\end{eqnarray*}
The last equivalence is obtained after the change of variable $y = \sqrt{d} x$ and $dy= dx\vert_{[1 : 0 : \dots : 0]}$ denotes the standard
Lebesgue measure of $\C^n$.

Now $P$ has been fixed once for all, so that $\int_{\C^n} \vert P(y) \vert^2 \exp (-\| y \|^2) dy $ is a constant. From now on I will
normalize $\sigma_P$ by setting
\begin{eqnarray}
\label{sigmaP}
\sigma_P &=& \sqrt{d}^n \frac{Q_d X_0^{d-k} }{\sqrt{\int_{\C^n} \vert P(y) \vert^2 \exp (-\| y \|^2) dy }}.
\end{eqnarray}
This section has $L^2$-norm one asymptotically, this $L^2$-norm being still concentrated in a ball of radius $R/\sqrt{d} $,
 but near the origin $ [1 : 0 : \dots : 0]$, its pointwise Fubini-Study norm is of the order $\sqrt{d}^n $.

Note that the same holds true for the section $\sqrt{\frac{(d+n) !}{n! d ! } } X_0^{d} $ above, which corresponds to $\sigma_P $ for $P=1$
(modulo the normalization of the volume) and this sheeds some light on the coefficients $\sqrt{\frac{(d+n) !}{n! \alpha_0 ! \dots \alpha_n !}} $ instead of 
$\sqrt{d \choose \alpha_0  \dots \alpha_n } $ in the orthonormal basis obtained in the above exercise and introduced before Theorem
\ref{theoGaWe}. 

\subsection{Probability of presence of $\Sigma$}

Recall that I did fix a closed hypersurface $\Sigma \subset B(0,R) \subset \R^n$ which does not need to be connected.
I then did construct, for every $x \in \R P^n$, a homogeneous polynomial $\sigma_P \in \R_d^{\text{hom}} [X_0 , \dots , X_n]$ such that $\sigma_P^{-1} (0) \cap B(x , R/\sqrt{d})$ contains a union of components $\widetilde{\Sigma}$ for which the pair $(B(x , R/\sqrt{d}) , \widetilde{\Sigma})$ is diffeomorphic to $(\R^n , \Sigma)$. I now claim much more. 

\begin{theo}[joint with Damien Gayet, \cite{GaWe4}]
\label{theoprob}
There exist $\tilde{c}_\Sigma > 0$ such that for every $x \in \R P^n$,
$$\liminf_{d \to + \infty} \mu \left\{ \sigma \in \R_d^{\text{hom}} [X_0 , \dots , X_n] \, \left| \, 
\begin{array}{l}
\sigma^{-1} (0) \cap B(x , R/\sqrt{d}) \supset \widetilde{\Sigma} \\
(B(x , R/\sqrt{d}) , \widetilde{\Sigma}) \cong (\R^n , \Sigma)
\end{array} \right. \right\} \geq \tilde{c}_\Sigma.$$
\end{theo}
Hence, not only there exists a polynomial $\sigma_P \in \R_d^{\text{hom}} [X_0 , \dots , X_n]$ with our desired properties, but moreover we
had a positive probability to find such, probability uniformely bounded from below by a positive constant. \\

{\bf Proof:}

{\bf First step:} Let me choose tubular neighborhoods $K$ and $U$ of $\Sigma$, $K$ being compact,  such that 
$\Sigma \subset K \subset U \subset B(0,R) $ and
\begin{enumerate}
\item $\vert P \vert_{U\setminus K} > \delta$, so that in particular $ P$ does not vanish in $U \setminus K$,
\item If $\vert P(y) \vert \leq \delta, y \in U$,  then $\vert d_y P \vert > \epsilon$,
\end{enumerate}
for some $\delta, \epsilon >0$.

Now, let me denote by $\Sigma_d ,  K_d ,  U_d$ the images of $\Sigma ,  K$ and $U$ under the homothety of rate $\frac{1}{\sqrt{d}}$,
so that  $\Sigma_d \subset K_d \subset U_d \subset B(0,R/\sqrt{d}) $.
I get likewise,
\begin{enumerate}
\item $\vert \sigma_P \vert_{U_d \setminus K_d} > \delta \sqrt{d}^n$ and
\item If $\vert \sigma_P (y) \vert \leq \delta \sqrt{d}^n, y \in U_d$,  then $\vert d_y \sigma_P \vert > \epsilon \sqrt{d}^{n+1}$,
\end{enumerate}
for some may be slightly different constants $\delta, \epsilon >0$.

This first step is called quantitative transversality. I knew that $0$ is a regular value of $\sigma_P $, but I am quantifying how much transversal to
the zero section $\sigma_P $ is. Such kind of quantitative transversality played a key role in the construction by S. K. Donaldson of symplectic divisors
in any closed symplectic manifold, see \cite{Don}.\\

{\bf Second step:}

\begin{prop}[joint with Damien Gayet, \cite{GaWe4}]
\label{propsup}
There exist $C_1 , C_2 > 0$ such that
$$\E (\sup_{B(0,R/\sqrt{d})} \vert \sigma \vert ) \leq C_1 \sqrt{d}^n \text{ and } \E (\sup_{B(0,R/\sqrt{d})} \vert d \sigma \vert) \leq C_2 \sqrt{d}^{n+1}. \; \square$$
\end{prop}
Let me skip the proof of this proposition, but point out however that
\begin{eqnarray*}
\E (\vert \sigma ([1 : 0 : \dots : 0]) \vert) &=& \int_{\R_d^{\text{hom}} [X_0 , \dots , X_n]} \vert \sigma ([1 : 0 : \dots : 0]) \vert d\mu (\sigma)\\
&=& \int_{\langle X_0^d \rangle} \vert \sigma ([1 : 0 : \dots : 0]) \vert d\mu (\sigma)\\
&\sim_{d \to + \infty} & \frac{\sqrt{d}^n}{\sqrt{n!}},
\end{eqnarray*}
if the volume of $\C P^n$ has been normalized to one.\\

{\bf Third step:} (from now on we follow an approach similar to the one used by Nazarov and Sodin in \cite{NazSod}).
Recall the following.

\begin{theo}[Markov's inequality]
Let $(\Omega , \mu)$ be a probabilty space and $f : \Omega \to \R^+$ be a random variable. Let $e = \E (f) = \int_\Omega f d\mu$
be its expectation. Then, for every $C >0$, $\mu \{ \omega \in \Omega \, \vert \, f(\omega) \geq C \} \leq e/C$.
\end{theo}

{\bf Proof:}

$$e = \int_\Omega f d\mu \geq \int_{\{ \omega \in \Omega \, \vert \, f(\omega) \geq C \} } f d\mu \geq C \mu \{ \omega \in \Omega \, \vert \, f(\omega) \geq C \} $$
$$\implies \mu \{ \omega \in \Omega \, \vert \, f(\omega) \geq C \} \leq e/C. \; \square$$\\

{\bf Application:} Since $\E (\sup_{B(0,R/\sqrt{d})} \vert \sigma \vert ) \leq C_1 \sqrt{d}^n$, from Markov's inequality we deduce that
$$\mu \{ \sigma \in \R_d^{\text{hom}} [X_0 , \dots , X_n] \, \vert \, \sup_{B(0,R/\sqrt{d})} \vert \sigma \vert \geq 4C_1 \sqrt{d}^n \} \leq \frac{1}{4}$$
and likewise
$$\mu \{ \sigma \in \R_d^{\text{hom}} [X_0 , \dots , X_n] \, \vert \, \sup_{B(0,R/\sqrt{d})} \vert d \sigma \vert \geq 4C_2 \sqrt{d}^{n+1} \} \leq \frac{1}{4},$$
so that
$$\mu \left\{ \sigma \in \R_d^{\text{hom}} [X_0 , \dots , X_n] \, \left| \, 
\begin{array}{l}
\sup_{B(0,R/\sqrt{d})} \vert \sigma \vert \leq 4C_1 \sqrt{d}^n\\
\sup_{B(0,R/\sqrt{d})} \vert d \sigma \vert \leq 4C_2 \sqrt{d}^{n+1}
\end{array}
\right.
\right\} \geq \frac{1}{2}.$$\\

{\bf Last step:}

Recall that I have to find a subset $E \subset \R_d^{\text{hom}} [X_0 , \dots , X_n]$ of measure uniformely bounded from below by some
positive constant, such that any polynomonial $\sigma$ in $E$ has the property that $\sigma^{-1} (0) \cap B(x , R/\sqrt{d}) \supset \widetilde{\Sigma}$ and $(B(x , R/\sqrt{d}) , \widetilde{\Sigma}) \cong (\R^n , \Sigma)$. This subset is going to be the set
$$
E_M = \left\{ a \sigma_P + \tau \, \left| \,
\begin{array}{l}
a \geq M\\
\tau \in \sigma_P^\perp\\
\sup_{B(0,R/\sqrt{d})} \vert \tau \vert \leq 4C_1 \sqrt{d}^n\\
\sup_{B(0,R/\sqrt{d})} \vert d \tau \vert \leq 4C_2 \sqrt{d}^{n+1}
\end{array}
\right.
\right\}, $$
with $M = \sup \{ \frac{4C_1}{\delta} , \frac{4C_2}{\epsilon} \}$.The measure of $E_M$ satisfies
\begin{eqnarray*}
\mu (E_M) &=& \big(\int_M^{+\infty} \exp (-t^2) \frac{dt}{\sqrt{\pi}} \big) \mu \left\{ \tau \in \sigma_P^\perp \, \left| \,\begin{array}{l}
\sup_{B(0,R/\sqrt{d})} \vert \tau \vert \leq 4C_1 \sqrt{d}^n\\
\sup_{B(0,R/\sqrt{d})} \vert d \tau \vert \leq 4C_2 \sqrt{d}^{n+1}
\end{array}
\right.
\right\}\\
&\geq& \frac{1}{2} \int_M^{+\infty} \exp (-t^2) \frac{dt}{\sqrt{\pi}}\\
&=& \tilde{c}_\Sigma >0.
\end{eqnarray*}
One checks indeed that 
$$\mu \left\{ \tau \in \sigma_P^\perp \, \left| \,\begin{array}{l}
\sup_{B(0,R/\sqrt{d})} \vert \tau \vert \leq 4C_1 \sqrt{d}^n\\
\sup_{B(0,R/\sqrt{d})} \vert d \tau \vert \leq 4C_2 \sqrt{d}^{n+1}
\end{array}
\right.
\right\} \geq  \frac{1}{2}$$
as
$$\mu \left\{ \sigma \in \R_d^{\text{hom}} [X_0 , \dots , X_n] \, \left| \, 
\begin{array}{l}
\sup_{B(0,R/\sqrt{d})} \vert \sigma \vert \leq 4C_1 \sqrt{d}^n\\
\sup_{B(0,R/\sqrt{d})} \vert d \sigma \vert \leq 4C_2 \sqrt{d}^{n+1}
\end{array}
\right.
\right\} \geq \frac{1}{2}.$$

Now, let $\sigma \in E_M$, $\sigma =  a \sigma_P + \tau$. For every $t \in [0,1]$, set
$\sigma_t =  a \sigma_P + t \tau$, so that $\sigma_0 =  \sigma_P$ and $\sigma_1 =  \sigma$.
Then, for every $t \in [0,1]$, $\sigma_t $ vanishes transversely in the open set $U_d$.
Indeed, let $x \in U_d$ and $t \in [0,1]$ such that $\sigma_t (x) = 0$. Then
\begin{eqnarray*}
\vert a \sigma_P  (x)  \vert & = & \vert  t \tau (x) \vert \leq 4C_1 \sqrt{d}^n\\
\implies  \vert  \sigma_P  (x)  \vert & \leq & \delta \sqrt{d}^n\\
\implies  \vert  d \sigma_P  (x)  \vert & > & \epsilon \sqrt{d}^{n+1},
\end{eqnarray*}
so that
\begin{eqnarray*}
 \vert d \sigma_t  \vert & = & \vert  a d \sigma_P + t d \tau (x) \vert \\
& \geq &  \vert  a d \sigma_P \vert   - \vert   d \tau (x) \vert \\
&>& 4C_2 \sqrt{d}^{n+1} - 4C_2 \sqrt{d}^{n+1}\\
&>& 0.
\end{eqnarray*}
We deduce that the smooth vanishing locus $\sigma_t^{-1} (0) \cap U_d = \Sigma_t$ remains trapped in the compact set $K_d$, thus cannot leave
$U_d$ and that it realizes an isotopy between $\Sigma = \sigma_P^{-1} (0) \cap U_d$ and $\sigma^{-1} (0) \cap U_d $. Hence the result. 
$\square$

\subsection{Proof of the lower estimates}
\label{subsecproof}

Let me address the following problem. \\

{\bf Question $3$:} Given a closed Riemannian manifold $(M, g)$ and $\epsilon >0$, how many disjoint balls of radius
$\epsilon$ can be packed in $M$?\\

In our case, $(M, g)$ is going to be the real projective space $\R P^n$ equipped with its Fubini-Study metric and $\epsilon $
is going to be $R/\sqrt{d}$. Note that if $\C P^n$ is the quotient of the unit sphere $S^{2n+1}$ under the action of the unit circle
$U_1 (\C)$ by isometries, giving rise to the Hopf fibration, $\R P^n$ is just the quotient of the sphere $S^{n+1}$ under the action
of the group $\{ \pm 1 \}$ of unit elements of $\R$. This antipodal action is also isometric for the round metric of $S^{n+1}$
and the Fubini-Study metric of $\R P^n$ is just the metric on the quotient $\R P^n = S^{n+1} / \{ \pm 1 \}$ induced by this round metric. 

\begin{prop}
\label{proppack}
Let $(M, g)$ be a closed Riemannian manifold of dimension $n$ and  $\epsilon >0$. Let $N_\epsilon$ be the maximal number of disjoint balls of radius $\epsilon$ that can be packed in $M$. Then,
$$\liminf_{\epsilon \to 0} (\epsilon^n N_\epsilon) \geq \frac{\text{Vol}_g (M)}{2^n \text{Vol}_{\text{eucl}} (B(0,1))},$$
where $\text{Vol}_g (M)$ denotes the total Riemannian volume of $M$ and $\text{Vol}_{\text{eucl}} (B(0,1))$ the Euclidean volume of the unit ball in $\R^n$.
\end{prop}
Note that it is of course not possible to fill more than the total volume of $M$ by disjoint balls, so that 
$$\limsup_{\epsilon \to 0} (\epsilon^n N_\epsilon) \leq \frac{\text{Vol}_g (M)}{\text{Vol}_{\text{eucl}} (B(0,1))},$$
but from Proposition \ref{proppack} we know that it is possible to fill a fraction of it. This packing problem is a classical one. For instance,
in the case of the Euclidean space $\R^n$, the question may be, given a box, can we fill its whole volume with apples. Of course not and actually even if the radius of the apples was converging to zero. The question then becomes what is the best way to fill the box in order to loose
the minimal amount of space, but we do not address this question. If instead of Euclidean balls, we just wanted to fill the manifolds with 
balls of a given volume, that is by disjoint images of embeddings of the Euclidean balls by diffeomorphisms which preserve the volume form,
then it would be possible to fill the whole volume, say asymptotically due to Moser's trick. Finally, a famous theorem of M. Gromov establishes
that it is not possible to fill the whole Fubini-Study volume of $\C P^2$ by packing two disjoint symplectic balls, see \cite{Gro}, \cite{McDPo}, meaning two disjoint embeddings of some ball of $\C^2$ into $\C P^2$ which preserve the symplectic form. \\

{\bf Proof:}

Let $\Lambda_\epsilon$ be a subset of points of $M$ with the property that for every $x \neq y \in \Lambda_\epsilon$, $d(x,y)> 2 \epsilon$
and that $\Lambda_\epsilon$ is maximal with respect to this property. Then, the balls centered at the points of $\Lambda_\epsilon$ and
of radius $\epsilon$ are disjoint to each other, so that $\# \Lambda_\epsilon \leq N_\epsilon$. But the balls centered at the points of $\Lambda_\epsilon$ and of radius $2 \epsilon$ cover $M$ since a point $y$ in the complement of these balls in $M$ could be added
to $\Lambda_\epsilon$ to get a strictly larger set with our desired property, contradicting the maximality of $\Lambda_\epsilon$.
Thus 
$$\text{Vol}_g (M) \leq \sum_{x \in \Lambda_\epsilon} \text{Vol}_g (B(x, 2 \epsilon)) \sim_{\epsilon \to 0} \# \Lambda_\epsilon \epsilon^n 2^n \text{Vol}_{\text{eucl}} (B(0,1)),$$
 so that $\liminf_{\epsilon \to 0} (\epsilon^n \# \Lambda_\epsilon) \geq \frac{\text{Vol}_g (M)}{2^n \text{Vol}_{\text{eucl}} (B(0,1))}$. $\square$\\

Let me now come back to the proof of the lower estimates in Theorem \ref{theoGaWe}. Let $\epsilon = R/\sqrt{d}$ and $\Lambda_\epsilon$
be a subset of $(\R P^n, g_{\text{FS}})$ maximal with the property that for every $x \neq y \in \Lambda_\epsilon$, $d(x,y)> 2 \epsilon$.
For every closed connected hypersurface $\Sigma$ of $\R^n$ and $P \in \R_d^{\text{hom}} [X_0 , \dots , X_n]$, let $N_\Sigma (V_P)$
be the number of connected components of $V_P = P^{-1} (0) \subset \R P^n$ which are diffeomorphic to $\Sigma$. For every $x \in \R P^n$,
we set $N_{\Sigma, x} (V_P)$ to be one if $V_P \cap B(x, \epsilon) \supset \widetilde{\Sigma}$ such that $(B(x, \epsilon) , \widetilde{\Sigma})
\cong (\R^n , \Sigma)$ and $N_{\Sigma, x} (V_P) = 0$ otherwise. We deduce in particular, $N_\Sigma (V_P) \geq \sum_{x \in \Lambda_\epsilon} 
N_{\Sigma, x} (V_P)$. Then,
\begin{eqnarray*}
\E (b_i) &=& \int_{\R_d^{\text{hom}} [X_0 , \dots , X_n]} b_i (V_P) d\mu(P)\\
&\geq & \int_{\R_d^{\text{hom}} [X_0 , \dots , X_n]} \big( \sum_{[\Sigma] \in {\cal H}_n} b_i (\Sigma) N_\Sigma (V_P) \big) d\mu(P)\\
&\geq & \sum_{[\Sigma] \in {\cal H}_n} b_i (\Sigma) \sum_{x \in \Lambda_\epsilon} \int_{\R_d^{\text{hom}} [X_0 , \dots , X_n]} N_{\Sigma, x} (V_P)
d\mu(P)\\
&\geq & \sum_{[\Sigma] \in {\cal H}_n} b_i (\Sigma) \tilde{c}_\Sigma \# \Lambda_\epsilon \text{ from Theorem \ref{theoprob}}.
\end{eqnarray*}
Note that $d$ being fixed, all the sums involved are finite. We then use Proposition \ref{proppack} to deduce
$$\liminf_{\epsilon \to 0} \big(\frac{\# \Lambda_\epsilon}{\sqrt{d}^n} \big) \geq \frac{\text{Vol}_{\text{FS}} (\R P^n)}{2^n \text{Vol}_{\text{eucl}} (B(0,R))}.$$
We finally set $$c_\Sigma = \frac{\tilde{c}_\Sigma}{2^n \text{Vol}_{\text{eucl}} (B(0,R))}$$ to get
$$\liminf_{d \to +\infty}  \big( \frac{\E (b_i) }{\sqrt{d}^n \text{Vol}_{\text{FS}} (\R P^n)}  \big) \geq  \sum_{[\Sigma] \in {\cal H}_n} c_\Sigma b_i (\Sigma)  = c_i^-. \; \; \; \;\square$$

Note that we have actually proved that the expected number $\E (N_\Sigma)$ of connected components of $V_P$ diffeomorphic
to $\Sigma$ satisfies
$$\liminf_{d \to +\infty}  \big( \frac{\E (N_\Sigma) }{\sqrt{d}^n \text{Vol}_{\text{FS}} (\R P^n)}  \big) \geq c_\Sigma.$$
We have even proved this lower estimate for a smaller quantity, the expected number of disjoint balls $B$ of $\R P^n$ that contain
a component $\widetilde{\Sigma}$ of $V_P$ for which the pair $(B, \widetilde{\Sigma})$ gets diffeomorphic to $(\R^n , \Sigma)$. \\

{\bf Acknowledgement:} 

The research leading to these results has received funding from the European Community's Seventh Framework Progamme 
([FP7/2007-2013] [FP7/2007-2011]) under grant agreement $\text{n}\textsuperscript{o}$ [258204].

\addcontentsline{toc}{part}{\hspace*{\indentation}Bibliography}


\bibliography{CIMPA}
\bibliographystyle{abbrv}

\vspace{0.7cm}
\noindent 
 Universit\'e de Lyon \\
CNRS UMR 5208 \\
Universit\'e Lyon 1 \\
Institut Camille Jordan \\
43 blvd. du 11 novembre 1918 \\
F-69622 Villeurbanne cedex\\
France\\
 welschinger@math.univ-lyon1.fr

\end{document}